\documentclass[12pt,amsfonts,amssymb]{article} 
 \newtheorem{theor}{Theorem}[section]

 \newtheorem{definition}[theor]{Definition}

\newtheorem{notation}[theor]{Notation}

\newtheorem{theorem}[theor]{Theorem}
\newtheorem{lemma}[theor]{Lemma}
\newtheorem{question}[theor]{Question}
\newtheorem{fact}[theor]{Fact}
\newtheorem{claim}[theor]{Claim}
\def\tempsi{\Delta_1}

\def\Dscr{\cal D}

\def\av{\rm avg }
\def\dom{\rm dom\  } 
 
\def\b1K{\mbox{\boldmath $K$}_{-1}}

\def\TP{\tau}

\def\acl{{\rm acl}} 
\def\dcl{{\rm dcl}}
\def\restrict{\mid}

\def\diag{{\rm diag\ }}

\def\tp{\rm tp }

\def\abar{{\bf a}} 
\def\xbar{{\bf x}} 
\def\bbar{{\bf b}} 
\def\cbar{{\bf c}} 
\def\dbar{{\bf d}} 
\def\ebar{{\bf e}} 
\def\fbar{{\bf f}} 
 
\def\ybar{{\bf y}} 

\def\ubar{{\bf u}} 
\def\vbar{{\bf v}}
\def\wbar{{\bf w}}
 \def\zbar{{\bf z}}

\title{Model Companions of $T_{\sigma}$ for stable $T$ } 
\author{
John T. Baldwin
\thanks{Partially supported by NSF grant
DMS 9803496.}
\\ Department of Mathematics, Statistics and  Computer
Science\\University of Illinois at Chicago
\and
Saharon Shelah\thanks{This is paper 759 in Shelah's bibliography. 
The authors thank
Mittag-Leffler Institute and the Binational Science Foundation for partial
support of this research.}\\ Department of Mathematics\\
Hebrew University of Jerusalem\\Rutgers University}

\begin{document} 

\maketitle

Let $T$ be a complete first order theory in a countable 
relational language $L$.  
We assume relation symbols have been added to make each formula equivalent
to a predicate. Adjoin a new unary function symbol $\sigma$ to obtain
the language $L_\sigma$; $T_{\sigma}$
is obtained by adding axioms asserting that $\sigma$ is an $L$-automorphism.

The modern study of the model companion of theories with an automorphism
has two aspects.  One line, stemming from Lascar \cite{Lascarba}, deals
with `generic' automorphisms of arbitrary structures.  A second, beginning
with Chatzidakis and Hrushovski \cite{ChatHrushovski} and questions of
Macintyre about the Frobenius automorphism is more concerned with specific
algebraic theories.  This paper is more in the first tradition: we find
general necessary and sufficient conditions for a stable first order theory
with automorphism to have a model companion. 

Kikyo  investigates the existence of model companions of $T_\sigma$ when
$T$ is unstable in  \cite{Kikyo}.
He also includes an argument of  Kudaibergenov showing  that
if $T$ is stable with  the finite cover property 
then $T_{\sigma}$ has no
model companion.  This argument was implicit in 
\cite{ChatPillay} and 
is a rediscovery of a theorem of Winkler \cite{Winkler} in 
the 70's.
We provide  necessary and sufficient conditions for
 $T_{\sigma}$ to have a model companion when $T$ is stable.
Namely, we introduce a new condition: $T_\sigma$ admits obstructions,
and show that $T_\sigma$ has a model companion iff and only if
$T_\sigma$ does not admit obstructions.  This condition is weakening 
of the finite cover property: if a stable theory $T$ has the finite cover
property then $T_{\sigma}$ admits obstructions.

Kikyo also proved that if $T$ is an unstable theory without the
independence property, $T_\sigma$ does not have a model-companion.
Kikyo and  Shelah  \cite{KikyoShelah} have improved this
by weakening the hypothesis to, $T$ has the strict 
order property.

For $p$ a type over $A$ and $\sigma$ an automorphism with $A \subset \dom p$,
$\sigma(p)$ denotes $\{\phi(\xbar,\sigma(\abar)):\phi(\xbar,\abar) \in p\}$.
References of the form  II.4.13. are to 
\cite{Shelahbook2nd}

\section{Example}

  In the following example we examine exactly why a particular
$T_\sigma$ does not have a model companion.  Eventually,
we will show that the obstruction illustrated here represents  the 
reason $T_\sigma$ (for stable $T$) can fail to have a model companion.
Let $L$ contain two binary relation symbols
$E$ and $R$ and unary predicates $P_i$ for $i<\omega$.  The theory $T$
asserts that $E$ is an equivalence relation with infinitely many infinite classes,
which are refined by $R$ into two-element classes.  Moreover, each $P_i$
holds only of elements from one $E$-class and  contains exactly 
one element
from each $R$ class of that $E$-class. 

 Now, $T_\sigma$ does not have a model companion.
To see this, let 
$\psi(x,y,z)$ be the formula: $E(x,z) \wedge E(y,z) 
\wedge R(x,y)\wedge x \neq y$.
Let $\Gamma$ be the $L_\sigma$-type in
 the variables $\{z\} \cup \langle  x_i y_i: i < 
\omega\rangle$ which asserts for each $i$,
 $\{\psi(x_i,y_i,z)\}$, 
 the sequence  $\langle  x_i y_i: i < 
\omega\rangle$ is $L$-indiscernible  
and for every  $\phi(x,\wbar)\in L(T)$:
$$(\forall \wbar
) \bigvee\{
\bigwedge_{i\in U} \phi(x_i,\wbar) 
\leftrightarrow \phi(y_i,\sigma(\wbar)): U \subseteq \lg(\wbar) + 3,  
|U|>  (\lg(\wbar) + 3)/2\}.$$
Thus if  $\langle  b_i c_i: i < \omega\rangle a$ realize $\Gamma$ in a model $M$,
 $$\sigma( \av (\langle  b_i: i < \omega\rangle/M))= 
\av (\langle  c_i:  i < \omega\rangle/M).$$
For any finite $\Delta \subset L(T)$, let $\chi_{\Delta,k}(\xbar,\ybar,z)$
 be the  conjunction of the $\Delta$-formulas satisfied by
$\langle  b_i c_i: i < k\rangle a$ where 
$\langle  b_i c_i: i < k\rangle a$ are an initial segment of 
a realization of  $\Gamma$.  Let $\theta_{\Delta,k}$ be the  sentence
$$(\forall x_0, \ldots x_{k-1}, y_0, \ldots y_{k-1},z) \chi_{\Delta,k}
( x_0, \ldots x_{k-1}, y _0, \ldots y_{k-1},z) \rightarrow $$
$$
(\exists x_0, y_0, x_1,y_1)[  \psi(x_0,y_0,z) 
\wedge \psi(x_1,y_1,z)
\wedge \sigma (x_1) = y_1].$$


We  now claim that if $T_{\sigma}$ has a model companion $T^*_\sigma$, then for some
$k$ and $\Delta$,
$$T^*_\sigma \vdash \theta_{\Delta,k}.$$
For this, let $M \models T^*_\sigma$ such that 
 $\langle b_i c_i:  i < k\rangle a$ satisfy $\Gamma$ in $M$.
Suppose $M \restrict L \prec N$ and $N$ is an $|M|^+$-saturated model of
$T$.  In $N$ we can find $b,c$ realizing the average of 
$\langle b_i:  i < \omega\rangle$ and  $\langle c_i:  i < \omega\rangle$
over $M$ respectively. 
Then $$\sigma( \av (\langle  b_i: i < \omega\rangle/M))= \av (\langle  c_i:  i < \omega\rangle/M)$$
and so there
 is an automorphism $\sigma^*$ of $N$ extending $\sigma$ and taking $b$ to $c$.  Since $(M,\sigma)$
is existentially closed ($T^*_\sigma$ is model complete), we can pull $b,c$
down to $M$.  By compactness, some finite subset $\Gamma_0$ of $\Gamma$ suffices
and letting $\Delta$ be the formulas mentioned in  $\Gamma_0$ and $k$
 the number of
$x_i,y_i$ appearing in  $\Gamma_0$ we have the claim.

But now we show that if $(M,\sigma)$ is any model of $T_\sigma$, then
for any finite $\Delta$ and any $k$,
 $(M,\sigma) \models\neg  \theta_{\Delta,k}$. For this, choose
$b_i,c_i$ for $i < k$ which are $E$-equivalent to each other and to an element 
$a$ in a class $P_j$ where $P_j$  does not occur  in $\Delta$
 and with $R(b_i,c_i)$ and $b_i \neq c_i$.
Then   $\bbar,\cbar, a$ satisfy $\chi_{\Delta,k}$ but there are no  $b_k,c_k$ and automorphism
$\sigma$ which  makes  $\theta_{\Delta,k}$ true.
For,  for each $j$,
$$T \vdash (\forall x,y,z) (\psi(x,y,z) \wedge P_j(z) \rightarrow
 [P_j(x) \leftrightarrow \neg P_j(y)]).$$

To put this situation in a more general framework, recall some
notation from \cite{Shelahbook2nd}.  $\Delta$ will note a finite
set of formulas: $\{\phi_i(\xbar,\ybar_i): \lg(\xbar) =m, i < |\Delta|\}$;
 $p$ is a 
$\Delta$-$m$-type over $A$ if $p$ is a set of formulas 
 $\phi_i(\xbar,\abar)$ where $\xbar = \langle x_1,\ldots x_{m-1}\rangle$
 (these specific variables) and $\abar$  from $A$ is substituted 
for $\ybar_i$.  Thus, if $A$
is finite there are only finitely many $\Delta$-$m$-types over $A$.

Now let $\Delta_1$ contain Boolean combinations 
 of $x =y, R(x,y), E(x,y)$.
Let $\Delta_2$ expand $\Delta_1$ by adding a finite number of the
$P_j(z)$ and let $\Delta_3$ contain $P_j(x)$ where $P_j$ does
 not occur in $\Delta_2$.

Now we have the following situation:  there exists a set $X=\{
b_0,b_1,c_0,c_1,a\}$,  $P_j(a)$ holds, all  5 are $E$-equivalent
and $R(b_i,c_i)$  for $i = 0, 1$ such that: 
\begin{enumerate}
\item  $\langle  \bbar_i \cbar_i: i \leq  2\rangle$ is $\Delta_2$-indiscernible
over $a$.
\item $\langle b_0 c_0,b_1c_1\rangle$ can be extended to an infinite set of
 indiscernibles 
 $\overline{ \bbar}\overline{ \cbar}$ which satisfy the following. 
\item $\psi(b_i,c_i,a)$.
\item  $\sigma(\av_{\Delta_2}(\overline{ \bbar}/M)=\av_{\Delta_2}(\overline{ \cbar}/M)$.
\item $\tp_{\Delta_1}(\bbar_2\cbar_2/X) \vdash 
\sigma( \tp_{\Delta_3}(\bbar_2/X))\neq 
\tp_{\Delta_3}(\cbar_2/X)$.
\end{enumerate}

We call a sequence like  $\langle  \bbar_i \cbar_i: i \leq  2\rangle\abar$
a  $(\sigma,\tempsi, \Delta_2, \Delta_3,n)$-obstruction
 over the emptyset.
In order to `finitize' the notions we will give below more technical formulations
of the last two conditions: we will have to discuss
obstructions over a finite set $A$.  In the example, the identity was the only
automorphism of the prime model.  We will have to introduce a third sequence 
$\bbar'$ to deal with arbitrary $\sigma$. But  this example demonstrates 
the key aspects of obstruction which are the second reason for 
$T_\sigma$ to  lack a model companion.

\section{Preliminaries}

In order to express the notions described in the example, we need several notions
from basic stability theory.  By working with finite sets of formulas
in a stable theory without the finite cover property we are able to
refine arguments about infinite sets of indiscernibles to arguments about
sufficiently long finite sequences.  Let  $\Delta$ be a finite set of formulas
which we
will assume to be closed under permutations of variables and negation; $\neg\neg\phi$
is identified with $\phi$.  Recall that an ordered sequence 
$ E=\langle \abar_i:i\in I\rangle$ is said to
be $(\Delta,p)$ indiscernible over $A$ if any two properly ordered
 $p$-element 
subsequences of $I$ realize the same $\Delta$-type over $A$. 
 It is $\Delta$-indiscernible
if it is  $(\Delta,p)$-indiscernible  for all $p$, or equivalently for 
all $p'$ with $p'$ at most
the maximum number of variables in a formula in $\Delta$. 
For any sequence $ E=\langle \abar_i:i\in I\rangle$ and $j \in I$ we
write $E_j$ for $ \langle \abar_i:i< j\rangle$.

We will rely on the following 
facts/definitions  from \cite{Shelahbook2nd}
to  introduce two crucial functions for this paper: $ F(\Delta,n)$ and
 $f(\Delta,n)$.
\begin{fact}
\label{fact1}
Recall that if $T$ is stable, then for every finite $\Delta \subset L(T)$ 
and $n < \omega$
there 
is a finite $\Delta' = F(\Delta,n)$ with $\Delta \subseteq \Delta' 
\subset L(T)$  and 
a $k^* = f(\Delta,n)$ such that 
\begin{enumerate}
\item A sequence  $\langle \ebar_i:i \in I\rangle$ of $n$-tuples
such that for $i<j$ and a finite set $A$, $\tp_{F(\Delta,n)}(e_j/E_iA) =
\tp_{F(\Delta,n)}(e_i/E_iA)$ and 
 $R(F(\Delta,n),2)(e_j/E_jA)=R(F(\Delta,n),2)(e_i/E_iA)$,
(whence,  
 $\tp_{F(\Delta,n)}(e_j/E_iA)$ is definable over $A$)
 is a sequence of $\Delta$-indiscernibles (II.2.17).
\item  For any  set of  $\Delta'$-indiscernibles over the empty set,
$E =\langle \ebar_i \colon i < k\rangle$ with $\lg (\ebar_i) = n$ and  $k \geq k^*$ for any 
 $\theta(\ubar,\vbar) \in \Delta$ and any $\dbar$ with $
\lg(\dbar) = 
\lg(\vbar) = m$  either  
$\{\ebar_i\colon \phi(\ebar_i,\dbar)\}$ or 
$\{\ebar_i\colon\neg \phi(\ebar_i,\dbar)\}$ has 
strictly less than $k^*/2$ elements. (II.4.13., II.2.20 ) 
\item This implies that, for appropriate choice of $k^*$, 
\begin{enumerate}
\item there is an  integer $m = m(\Delta,n)\geq n$ such that
 for any  set of  $\Delta'$-indiscernibles
$\langle \ebar_i \colon i < k\rangle$ over  $A$ with $\lg (\ebar_i) = n$ and  $k \geq k^*$ and any
$\abar$ with $\lg(\abar)\leq m$
there is a $U \subseteq k$ with $|U|< k^*/2$ such that
$\langle \ebar_i \colon i \in  k - U\rangle$ 
is $\Delta$-indiscernible over $A\abar$;

\item moreover
 if $k \geq k^*$, for any set $A$, 
$\av_{\Delta}(\langle \ebar_i \colon i < k\rangle/A)$
is well-defined.  Namely,
 $\av_{\Delta}(\langle \ebar_i \colon i < k\rangle/A)
= $
$$\{\phi(\xbar,\abar): |\{\ebar_i:i<k, \phi(\ebar_i,\abar)| \geq \frac{k^*}{10},
\abar \in A, \phi(\xbar,\ybar) \in \Delta\}.$$
\end{enumerate}
\end{enumerate}
\end{fact}
In a), $m$ is the least $k \geq n$ such all $\phi\in \Delta$ have at
most $k$ free variables.
But $\av_{\Delta}(\langle \ebar_i \colon i < k\rangle/A)$ need not
be consistent.  (Let $E$ be all the members of one finite class in
the standard fcp example and let $A = E$.)  However,

\begin{fact}
\label{fact2}
If, in addition to Fact~\ref{fact1},  $T$ does not have the finite
 cover property, we can further demand

\begin{enumerate}
\item If
$E =\langle \ebar_i \colon i < k^*\rangle$ is a
 set of  $\Delta'$-indiscernibles over the empty set,
for any $A$, $\av_{\Delta}(E/A) $ is a {\em consistent} complete $\Delta$-type over
$A$.

\item  Moreover $k^*$ can be chosen so that any set of $\Delta'$-indiscernibles (of $n$-tuples)
with length at least $k^*$ can be extended to one of 
 infinite length (II.4.6).

\item For any  pair  of  $F(F(\Delta, n))$-indiscernible sequences
$E^1=\langle \ebar^1_i \colon i < k\rangle$ and 
$E^2=\langle \ebar^2_i \colon i < k\rangle$ over $\abar$  with $\lg (\ebar_i) = n$
and $k \geq k^*$ such that
$$ \av_{F(\Delta})(E^1/\abar E^1E^2)
=\av_{F(\Delta})(E^2/\abar E^1E^2),$$
 there exist ${\bf J} = \langle \ebar_j:k < j < \omega\rangle$ such that both
$E^1J$ and
$E^2J$ are $F(\Delta)$-indiscernible over $\abar$.
We express the displayed condition on 
$\overline{\ebar}^1,\overline{\ebar}^2$ by the
formula:
$\lambda_{\Delta}(\overline{\ebar}^1,\overline{\ebar}^2,\abar)$.

\item If $E_1$ and $E_2$ contained in a model $M$
are $F(\Delta,n)$-indiscernible over $\abar\in M$ and each have
length at least $k^*$, there is a formula 
$\lambda_\Delta(\xbar^1, \xbar^2,\zbar)$ such that
$M \models \lambda_\Delta(\ebar^1, \ebar^2,\abar)$ if and only if
$\av_\Delta(\ebar^1/M) = \av_\Delta(\ebar_2/M)$.
\end{enumerate}
\end{fact}

Proof.  For 1, make sure that $k^*$ is large enough that
every $\Delta$-type which is $k^*$-consistent is 
consistent (II.4.4 3)).
Now 3) follows by extending the common $F(\Delta)$-average of 
$E^1$ and $E^2$ over $\abar E^1E^2$ by 2).  Finally,
condition 4 holds by adapting the argument for III.1.8
from the set of all $L$-formulas to $\Delta$; 
$\lambda_\Delta$
is the formula from 3).

Note that both $F$ and $f$ can be chosen  increasing in $\Delta$ and $n$.

\section{Obstructions}
In this section we introduce the main new notion of this paper:
obstruction.

We are concerned with a formula $\psi(\xbar,\ybar,\zbar)$ where
$\lg(\xbar) = \lg(\ybar) =n$ and $\lg(\zbar) =m$.
We will apply Facts~\ref{fact1} and \ref{fact1}
with $\ebar_i = \bbar_i\bbar'_i\cbar_i$ 
where each of $\bbar_i$, $\bbar'_i$, and $\cbar_i$ has length $n$.
Thus, our exposition will depend on functions
 $F(\Delta, 3n)$, $f(\Delta, 3n)$.
In several cases, we apply Fact~\ref{fact1} with
$\phi(\ubar_1,\ubar_2,\ubar_3,\vbar)$ as 
$\theta(\ubar_2,\vbar) \leftrightarrow \theta(\ubar_3,\vbar)$ for
various $\theta$.  

The following notation is crucial to state the definition.

\begin{notation}
\label{taudef}
If $\overline {\dbar} = \langle \dbar_i: i< r\rangle$ is a sequence of 
$3n$-tuples, which is $\Delta$-indiscernible over
a finite sequence $\fbar$, and
  $r \geq k^* = f(\Delta,3n)$, then 
$\TP_{\Delta}( \zbar, \overline {\dbar}\fbar)$
is the formula with free variable $\zbar$ and 
parameters $\overline {\dbar}\fbar$ which asserts that
 there is a subsequence $\overline {\dbar}'$ of $\overline {\dbar}$
with length $f(\Delta,3n)$ 
  so that  $\overline {\dbar}'\zbar$ forms a 
sequence of $\Delta$-indiscernibles over $\fbar$. 
\end{notation}

The following result follows easily from  Fact~\ref{fact2} 3) and the definition of
$f(\Delta,n)$. 


\begin{lemma}  Let $p \geq k = f(\Delta,n)$ and suppose 
$\langle \ebar_i: i < p\rangle$ is 
a sequence of $3n$-tuples which is $F(\Delta,3n)$-
indiscernible over $\fbar$.  If $S_1, S_2 \subseteq p$ with
$|S_1|, |S_2| \geq k^*/2$ then
 $\TP_{\Delta}( \zbar, \overline {\dbar}|S_1\fbar)$ and 
 $\TP_{\Delta}( \zbar,\overline {\dbar}|S_2\fbar)$ are equivalent formulas.
\end{lemma}

Now we come to the main notion.
Intuitively, $\langle  \bbar_i \bbar'_i\cbar_i\colon i \leq k\rangle\abar$ 
is  a $(\Delta_1, \Delta_2, \Delta_3,n)$ 
obstruction over $A$
if $\langle  \bbar_i \bbar'_i\cbar_i\colon i \leq k\rangle$
is an indefinitely extendible sequence 
of $\Delta_2$ indiscernibles over $\abar$ such that
the $\bbar$'s, $\bbar'_i$'s
and
$\cbar$'s each have length $n$ and  the $\Delta_2$-average of the
$\bbar'_i$
and
$\cbar$ is the same (over any set)
but any realizations of the $\Delta_1$-type of the $\bbar'_i$ and
the $\Delta_1$-type of the $\cbar_i$ over $\abar$ and the sequence
have different $\Delta_3$-types over $A$.
More formally, we define:

\begin{definition}
\label{obstruction}
 For finite $\Delta_1 \subseteq \Delta_2 \subseteq L(T)$ and 
$\Delta_3 \subseteq L(T)$,
finite $\abar\subseteq A \subset M \models T$
with  $\lg(\abar) \leq m(\Delta,n)$ (as in Fact~\ref{fact1}),
 $\sigma$ an automorphism of $M$,
and a natural number  $n$,
 $\langle  \bbar_i \sigma(\bbar_i) \cbar_i: i \leq k\rangle\abar$ is  a $(\sigma, \Delta_1, \Delta_2, \Delta_3,n)$ 
{\em obstruction over} $A$
 if the following conditions hold.
\begin{enumerate}
\item  $\langle  \bbar_i\sigma(\bbar_i) \cbar_i: i \leq  k\rangle$ is $F(\Delta_2, 3n)$-indiscernible
over $\abar$.
\item $k\geq  f(\Delta_2,3n)$. 
\item $\av_{\Delta_2}(\overline{\ebar}^1/M) = \av_{\Delta_2}  (\overline{\ebar}^2/M)$
where $ (\overline{\ebar}^1 =\langle \sigma(\bbar_i): i< k\rangle$ and $\overline{\ebar}^2= \langle     {\cbar_i}: i< k\rangle$.

\item 
Writing  $\TP_{\Delta_1}(\langle \bbar_i\sigma(\bbar_i)\cbar_i):i<k\rangle 
\abar)$ with
free variables $\xbar, \xbar', \ybar$, we have
$$M \models (\forall \xbar,\xbar',\ybar)[\TP_{\Delta_1}( 
\xbar,\xbar',\ybar,\langle \bbar_i\sigma(\bbar_i)\cbar_i:i<k\rangle \abar) 
\rightarrow$$
$$
\bigvee\{\phi(\xbar',\fbar)\wedge \neg \phi(\ybar,\fbar): 
 \fbar \in A\cup \bigcup_{i<k}\langle 
\bbar_i\sigma(\bbar_i)\cbar_i\rangle \cup \abar, \phi \in \Delta_3\}]$$




\end{enumerate}

\end{definition}
By Fact~\ref{fact2}, Condition 3) is 
 expressed by a formula of $\ebar^1, \ebar^2$
and $\abar$. 
Crucially,  the hypothesis of the fourth condition in Definition~\ref{obstruction}
is an $L$-formula with parameters  $\langle \bbar_i\sigma(\bbar_i)\cbar_i):i<k\rangle \abar$;
the conclusion is an $L$-formula with parameters from $A$ as well.
The disjunction in the conclusion of condition 4) is nonempty since
each $\bbar_i$, $\sigma(\bbar_i)$ is in the domain if $\Delta_3$ is
nontrivial.  $\Delta_1$ and
$\Delta_2$ have $3n$ type-variables; $\Delta_3$ has $n$ type-variables.




\begin{fact}
\label{monot} Note that if  $\langle  \bbar_i \cbar_i: i \leq k\rangle \abar$ is  a $(\sigma,\Delta_1, \Delta_2, \Delta_3,n)$ 
obstruction over $A$ and $\Delta_1 \subseteq
\Delta'_2 \subseteq \Delta_2$, then 
 $\langle  \bbar_i \cbar_i: i \leq k\rangle \abar$ is  a $(\sigma,\Delta_1, \Delta'_2, \Delta_3,n)$ 
obstruction over $A$.  Further, if $\langle  \bbar_i \cbar_i: i \leq k\rangle \abar$ is  a $(\sigma,\Delta_1, \Delta_2, \Delta_3,n)$ 
obstruction over $A$ and $A\subseteq A'$, where $A'$ is finite,  
$\Delta_1 \subseteq \Delta'_2 \subseteq \Delta_2$, then 
 $\langle  \bbar_i \cbar_i: i \leq k\rangle \abar$ is  a $(\sigma,\Delta_1, \Delta'_2, \Delta_3,n)$ 
obstruction over $A'$.
\end{fact}

\begin{definition} 
\begin{enumerate}
\item We say $(M,\sigma) \models T_\sigma$ 
{\em has no $\sigma$-obstructions } when there is a function $G(\Delta_1,n)$ with
$F(\Delta_1,3n) \subseteq  G(\Delta_1,n)$  such that  if $\Delta_1$ is a finite subset of $L(T)$
  and  $ G(\Delta_1,n)$ is contained in the finite
$\Delta_3 \subset L(T)$, then for every finite subset $A$ of $M$, there is
no   $(\sigma, \Delta_1, G(\Delta_1,n), \Delta_3,n)$ 
obstruction over $A$.
\item We say $T_\sigma$ {\em has no  $\sigma$-obstructions} when
 there is a function $G(\Delta_1,n)$ 
(which does not depend on $(M,\sigma)$) such that for 
each $(M,\sigma)\models T_\sigma$,
 if $\Delta_1$ is a finite subset of $L(T)$, 
 $A$ is finite subset of $M$, and 
  $\Delta_3 $ is a finite subset of $L(T)$, there is no $(\sigma, \Delta_1, G(\Delta_1,n), \Delta_3,n)$ 
obstruction over $A$.
\end{enumerate}
\end{definition}

\begin{definition}  A simple obstruction is an obstruction where the
automorphism $\sigma$ is the identity.  The notions of a theory or model
having a simple obstruction are the obvious modifications of the previous
definition.
\end{definition}

\begin{lemma}  $T_\sigma$ has  obstructions if and only if
$T$ has simple obstructions.
\end{lemma}

Proof.  Suppose $T$ has obstructions; we must find simple obstructions.
 So, suppose for some
 $ \Delta_1$, and $n$, and for 
 every finite $\Delta_2 \supseteq F(\Delta_1,3n)$, there is 
a finite $\Delta_3 $ and
a tuple  $(M^{\Delta_2},\sigma^{\Delta_2},
 A^{\Delta_2},
k^{\Delta_2})$ such that:
$(M^{\Delta_2},\sigma^{\Delta_2}) \models T_\sigma$,  
$A^{\Delta_2}$ is a finite subset
of $M^{\Delta_2}$ and  $\bbar^{\Delta_2},  \sigma(\bbar^{\Delta_2}), 
\cbar^{\Delta_2},\abar^{\Delta_2}$ contained in 
 $M^{\Delta_2}$ are a  $(\sigma^{\Delta_2},\Delta_1, \Delta_2, \Delta_3,n)$
 obstruction of length $k^{\Delta_2}$ over 
 $A^{\Delta_2}$. Without loss of generality
 $\lg(\abar) = m = m(\Delta_1,3n)$ and 
we can write $\Delta_3 = \Delta_3(\Delta_2)$.
Now, define a family of simple obstruction by replacing each component
of the given sequence of obstructions by an appropriate object with left
prefix sim.  
$$^{{\rm sim}}A^{\Delta_2} = A^{\Delta_2} \cup
 \{\bbar^{\Delta_2},  \sigma(\bbar^{\Delta_2}), 
\cbar^{\Delta_2}\}$$
 $$^{{\rm sim}}(\bbar^{\Delta_2})=^{{\rm sim}}(\bbar'^{\Delta_2})
= \sigma(\bbar^{\Delta_2})$$
$$^{{\rm sim}}\cbar^{\Delta_2}=\cbar^{\Delta_2}.$$
We use the same sequence of formulas for the $\Delta_2$ and
$\Delta_3(\Delta_2)$.
It is routine to check that we now have an obstruction with respect to
the identity.

\begin{lemma}
\label{fcpobs} If $T$ is a stable theory with the finite cover property
then $T$ has obstructions.
\end{lemma}

Proof.  By II.4.1.14 of \cite{Shelahbook2nd}, there is a formula $E(\xbar,
\ybar,\zbar)$ such that for each $\dbar$, $E(\xbar,
\ybar,\zbar)$ is an equivalence relation and for each $n$ there is a 
$\dbar_n$ such that $E(\xbar,
\ybar,\dbar_n)$ has finitely many classes but more than $n$.  Let $\Delta_1$
be $\{E(\xbar,\ybar,\zbar), \neg E(\xbar,
\ybar,\zbar)\}$ and consider any $\Delta_2$. Fix $\lg(\xbar)=\lg(\ybar)=r$.
 There are arbitrarily large
sequences $\bbar_n =\langle \bbar^n_j:j<n\rangle,
 \cbar_n=\langle \cbar^n_j:j<n\rangle$
 such that for some $\dbar_{n_1}, \dbar_{n_2}$, 
$\bbar_n$
is a set of representatives for distinct classes of $E(\xbar,
\ybar,\dbar^n_1)$ while $\cbar_n$
is a  set of representatives of distinct classes for $E(\xbar,\ybar,\dbar^n_2)$ and $n_1 < n_2$.  So by compactness and Ramsey, 
for any $\Delta_2$ we  can find such $\bbar_k, \cbar_k$ where $k
= f(\Delta_2, 2r)$ and  $\bbar_k,\cbar_k$ are a sequence of length $k$
of $\Delta_2$ indiscernibles.  Now, if $\Delta_3$ contains formulas which
express that the number of equivalence classes of $E(\xbar,
\ybar,\zbar)$ is greater than $n_1$ and $A$ contains 
$\dbar_{n_1}, \dbar_{n_2}$, we have an
$({\rm identity}, \Delta_1, \Delta_2, \Delta_3,r)$ 
obstruction over $A$.  (Let $\bbar =\bbar' = \bbar_{k}$, $\cbar =\cbar_k$.)

\section{Model Companions of $T_{\sigma}$}

In this section we establish necessary and sufficient conditions on stable
$T$ for $T_{\sigma}$ to have a model companion.  First, we notice when
the model companion, if it exists, is complete.

 Note that
  $\acl(\emptyset) = \dcl(\emptyset)$ in $C^{eq}_T$ means every finite equivalence
relation $E(\xbar,\ybar)$ of $T$ is defined by a finite conjunction: $\bigwedge_{i< n} \phi_i(\xbar) \leftrightarrow 
\phi_i(\ybar).$

\begin{fact}
\label{fact0} 
\begin{enumerate}
\item  If $T$ is stable, $T_{\sigma}$ has the amalgamation property.
\item If, in addition,    $\acl(\emptyset) = \dcl(\emptyset)$ in $C^{eq}_T$
then $T_\sigma$ has the joint embedding property. 
\end{enumerate}
\end{fact}

Proof. The first part of this Lemma was  proved by 
 (Theorem 3.3 of  \cite{Lascarba})
using the definability of types. For the second part, the hypothesis implies that types over
 the empty set are stationary
and the result follows by similar arguments.

\begin{lemma} Suppose $T$ is stable and 
$ T_\sigma$ has a model companion $T^*_{\sigma}$.
\begin{enumerate}
\item Then $T^*_{\sigma}$ is complete if and only if
$\acl(\emptyset) = \dcl(\emptyset)$ in $C^{eq}_T$.
\item If $(M,\sigma) \models T_\sigma$ 
then the union of the complete diagram of $M$ (in $L$)
with the diagram of $(M,\sigma)$ and $T^*_{\sigma}$ is complete.
\end{enumerate}
\end{lemma}

Proof. 1)
We have just seen that if $\acl(\emptyset) = \dcl(\emptyset)$ in $C^{eq}_T$,
then $T_\sigma$ has the joint embedding property; this implies in general
that the model companion is complete.  If  $\acl(\emptyset) \neq  \dcl(\emptyset)$ in $C^{eq}_T$,
let $E(\xbar,\ybar)$ be a finite equivalence relation witnessing 
 $\acl(\emptyset) \neq \dcl(\emptyset)$.  Because $E$ is a finite equivalence relation, 
$$T_1 = T_{\sigma} \cup \{(\forall \xbar) E(\xbar,\sigma(\xbar))\}$$
is a consistent extension of $T_\sigma$.  But since 
$$T_\sigma \cup \{\neg E(\xbar,\ybar)\} \cup \{\phi(\xbar) \leftrightarrow \phi(\ybar)\colon \phi \in L(T)\}
$$
is consistent so is  
$$T_2 = T_{\sigma} \cup \{(\exists  \xbar)\neg E(\xbar,\sigma(\xbar))\}.$$

But $T_1$ and $T_2$ are contradictory, so $T^*_\sigma$ is not complete.

2) Since we have joint embedding (from amalgamation over any model) the result follows as in 
Fact~\ref{fact0}.

We now prove the equivalence of three conditions: 
the first is a condition on a pair of models.  The second
is given by an infinite set of $L_\sigma$ sentences (take the union
over all finite $\Delta_2$) and the average requires names for all
elements of $M$. The third is expressed by a single first order 
sentence in $L_\sigma$. 
 The equivalence of the 
first and third suffices (Theorem~\ref{final}) to show the existence of a 
model companion.  In fact 1 implies 2 implies 3 requires only stability;
the nfcp is used to prove 3 implies 1.

\begin{lemma} 
\label{maindeal}
Suppose $T$ is stable without the fcp. 
Let $(M,\sigma) \models T_\sigma$, $\abar \in M$ 
and suppose that $(M,\sigma)$ has no $\sigma$-obstructions. Fix $\psi(\xbar,\ybar,\zbar)$ with
$\lg(\xbar) = \lg(\ybar)= n$ and $\lg(\zbar) = \lg (\abar)= m$. 

The following three assertions are equivalent:
\begin{enumerate}
\item
There exists $(N,\sigma)$,   $(M,\sigma) \subseteq (N,\sigma) \models T_\sigma$
and $$N \models (\exists \xbar \ybar)[\psi(\xbar,\ybar,\abar) \wedge \sigma(\xbar) =\ybar].$$
\item Fix  $\Delta_1 =\{\psi(\xbar,\ybar,\zbar)\}$
and without loss of generality $\lg(\zbar) \leq m(\Delta_1,n)$.
 For $k \geq 5 \cdot f(\Delta_1, 3n)$
any finite  $\Delta_2 \supseteq  F(\Delta_1,3n)$ (Fact~\ref{fact1}),
there are $\bbar_i,\sigma(\bbar_i)\cbar_i \in\, ^{3n}M$ for $i<k$ such that
\begin{enumerate}
\item $\langle  \bbar_i \sigma(\bbar_i) \cbar_i: i < k\rangle$ is
 $F(\Delta_2,3n)$-indiscernible over $\abar $,  
\item for each $i< k$, $\psi(\bbar_i, \cbar_i,\abar)$ holds, 
\item For every $\dbar \in\, ^mM$ and $\phi(\ubar,\vbar) \in \Delta_2$ we have
$$|\{i<k:\phi(\sigma(\bbar_i),\dbar) \leftrightarrow \phi(\cbar_i,\dbar)\}| \geq f(\Delta_2,3n)/2.$$
\end{enumerate}
\item Let $\Delta_2 = G(\Delta_1,n)$. Then are there are $\bbar_i,\sigma(\bbar_i)\cbar_i 
\in\, ^{3n}M$
 for $i<k =5 \cdot f(\Delta_2,3n)$ such that :
\begin{enumerate}
\item $\langle  \bbar_i \sigma(\bbar_i) 
\cbar_i: i < k\rangle$ 
is $F(G(\Delta_1,n), 3n)$-indiscernible over $\abar$,
\item for each $i< k$, $\psi(\bbar_i, \cbar_i,\abar)$,
\item $\lambda_{\Delta_2}(\langle \sigma(\bbar_i):i<k\rangle,
\langle \cbar_i:i<k\rangle,\abar)$.
\end{enumerate}

\end{enumerate}
\end{lemma}

Proof. 
 First we show 1) implies 2). Fix $\bbar,\cbar \in N$ with
 $N \models \psi(\bbar,\cbar,\abar) \wedge \sigma(\bbar) =\cbar$.
For  $\Delta_2$, let $\Delta^+_2 = F(F(\Delta_2,3n),3n)$.
 For each $\Delta_2$, 
choose a finite $p \subseteq \tp_{L(T)}(\bbar,\cbar/M)$
 with the same $(\Delta^+_2,2)$ rank  as $\tp_{\Delta_2}(\bbar,\cbar/M)$ (so
  $\tp_{\Delta_2}(\bbar,\cbar/M)$ is definable over $\dom p$). 
Now inductively construct an 
$F(\Delta_2,3n)$-indiscernible sequence  (by Fact~\ref{fact1} 1))  $\langle \bbar_i,\cbar_i: i< \omega\rangle$
by choosing $\bbar_i,\cbar_i$ in $M$ realizing the restriction of 
   $\tp_{\Delta^+_2}(\bbar,\cbar/M)$ to $\dom p$ along with the points 
already chosen. Let $\bbar'_i = \sigma(\bbar_i)$.
By  Ramsey's Theorem for some infinite 
 $U \subseteq \omega$, 
$\langle \bbar_i,\bbar'_i,\cbar_i: i \in U\rangle$ is 
$F(\Delta_2,3n)$-indiscernible
over $\abar$; renumbering let $U = \omega$.  Now conditions  a) and b)  of assertion 2) are clear.
For clause c), 
$$\av_{\Delta_2}(\langle \cbar_i: i<\omega\rangle/M) = 
\tp_{\Delta_2}(\cbar,M)$$ 
$$ =\sigma(\tp_{\Delta_2}(\bbar,M))
=\sigma(\av_{\Delta_2}(\langle \bbar_i: i<\omega\rangle/M))$$
since $\sigma(\bbar) = \cbar$.
So, for  each $\phi \in \Delta_2$  and each $\dbar \in M$ of 
appropriate length,
$$\phi(\xbar,\dbar) \in \av_{\Delta_2}(\langle \cbar_i:
i<\omega\rangle/M)$$
if and only if 
$$\phi(\xbar,\sigma^{-1}(\dbar)) \in \av_{\Delta_2}(\langle \bbar_i:
i<\omega\rangle/M).$$
So for some $S_1, S_2 \subset \omega$ with $|S_1|, |S_2| < f(\Delta_2,3n)/2$,
we have for all  $i \in \omega -(S_1 \cup S_2)$,
$\phi(\cbar_i,\dbar)$ if and only if $\phi(\bbar_i,\sigma^{-1}(\dbar))$.
Since $\sigma$ is an automorphism of $M$ this implies  for
 $i \in \omega -(S_1 \cup S_2)$,
$\phi(\cbar_i,\dbar)$ if and only if $\phi(\sigma(\bbar_i),\dbar)$
which gives condition c) by using the first $k$ elements of 
 $\langle \bbar_i, \sigma(\bbar_i)\cbar_i: i< \omega - S_1 \cup S_2\rangle$.

3) is a special case of 2). To see this, note that 3c) is easily implied
by the form analogous to 2c): For every $m  \leq m(\Delta_1,n)$ and   $\dbar \in \,^mM$ and 
$\phi(\ubar,\vbar) \in G(\Delta_1,n)$  we have
$$|\{i<k:\phi(\sigma(\bbar_i),\dbar) \leftrightarrow 
\phi(\cbar_i,\dbar)\}| \geq f(\Delta_2,3 n)/2.$$ 
If $T$ does not have f.c.p. the converse holds and we use that fact implicitly
in the following argument.
It remains only to show that 3) implies 1)
with $\Delta_1 = \{\psi\}$ and $\Delta_2 = G(\Delta_1,n)$.
Without loss of generality we may assume $N$ is $\aleph_1$-saturated.
We claim the type $$\Gamma=\{\psi(\xbar,\ybar,\abar)\} 
\cup \{\phi(\xbar,\dbar) \leftrightarrow
\phi(\ybar,\sigma(\dbar)): \dbar \in M, \phi\in L(T)\} \cup 
\diag(M)$$ is consistent.  This clearly suffices.

Let $k =f(\Delta_2,3n)$.  Suppose $\langle 
 \bbar_i \sigma( \bbar_i) \cbar_i: i <k\rangle\abar$
 satisfy 3).  
Let $\Gamma_0$ be a finite subset of $\Gamma$ and suppose only 
formulas from the
finite set  $\Delta_3 $ 
and only parameters from the finite set $A$  appear in $\Gamma_0$.
Write $\bbar'_i$ for $\sigma(\bbar_i)$.

Now $\langle  \bbar_i  \bbar'_i \cbar_i: i \leq f(\Delta_2,3n)\rangle\abar$
easily satisfy the first two conditions of Definition~\ref{obstruction}
for being a $(\Delta_1,\Delta_2,\Delta_3,n)$-obstruction over $A$ and,
in view of Fact~\ref{fact2} 3), 4), the 
third is given by condition 3c).
Since there is no obstruction the 4th condition must fail.  
So there exist $\bbar^*, (\bbar^*)', \cbar^*$ so that 
$$M \models \TP_{\Delta_1}(\bbar^*,(\bbar*)',\cbar^*,
\langle \bbar_i,\bbar'_i,\cbar_i: i < k\rangle\abar).$$ and
$\tp_{\Delta_3}((\bbar^*)'/A) = \tp_{\Delta_3}(\cbar^*/A)$ so
$\Gamma_0$ is satisfiable.

As we'll note in Theorem~\ref{final}, we have established a sufficient
condition for $T_\sigma$ to have a model companion.  The next argument
shows it is also necessary.

\begin{lemma}
\label{obsneeded} Suppose $T$ is stable; if $T$
has an obstruction then $T_{\sigma}$ does not have a model companion.

More precisely, suppose for some
 $ \Delta_1$, and $n$, and for 
 every finite $\Delta_2 \supseteq F(\Delta_1,3n)$, there is a 
a finite $\Delta_3 $ and
a tuple  $(M^{\Delta_2},\sigma^{\Delta_2},
 A^{\Delta_2},
k^{\Delta_2})$ such that:
$(M^{\Delta_2},\sigma^{\Delta_2}) \models T_\sigma$,  
$A^{\Delta_2}$ is a finite subset
of $M^{\Delta_2}$, $\bbar^{\Delta_2},  \sigma(\bbar^{\Delta_2}), 
\cbar^{\Delta_2},\abar^{\Delta_2}$ contained in 
 $M^{\Delta_2}$ are a  $(\sigma^{\Delta_2},\Delta_1, \Delta_2, \Delta_3,n)$
 obstruction of length $k^{\Delta_2}$ over 
 $A^{\Delta_2}$.
Without loss of generality $\lg(\abar) = m = m(\Delta_1,3n)$ and 
we can write $\Delta_3 = \Delta_3(\Delta_2)$.

Then the collection ${\bf K}_\sigma$ of existentially closed models of $T_\sigma$
 is not an elementary class.

\end{lemma}

Proof.  We may assume $T$ does not have f.c.p., since if it does we
know by Winkler and Kudaibergerov that $T_\sigma$ does not have 
a model companion.
By the usual coding we may assume  
$\Delta_1 = \{\psi(\xbar,\ybar,\zbar)\}$ with 
$\lg (\xbar)=\lg (\ybar)= n$, $\lg(\zbar) =m$, 
$k = f(\Delta_1,3n)$ .  Without loss of generality each $M^{\Delta_2}$ is existentially closed.
Let $\Dscr$ be a nonprincipal ultrafilter on
 $Y = \{\Delta_2: F(\Delta_1,3n)  \subseteq \Delta_2 \subset_\omega L(T)\}$
such that for any $\Delta \in Y$ the family of supersets in $Y$ of $\Delta$ is in $\Dscr$.
Expand the language to $L^+$ by adding a  unary function symbol $\sigma$ 
and a $3n$-ary relation symbol $Q$
and constants $\abar$.  Expand
each of the  $M^{\Delta_2}$ to an  $L^+$-structure $N_{\Delta_2}$ by
interpreting  
$\abar$ as $\abar^{\Delta_2}$
 $\sigma$ as $\sigma^{\Delta_2}$ and
$Q$ as the set  $\{\bbar_i^{\Delta_2},\sigma(\bbar_i^{\Delta_2}) \cbar_i^{\Delta_2}: i< k^{\Delta_2}\}$
of $3n$-tuples.
Let $(N^*,\sigma^*,Q^*,\abar^*)$ be the ultraproduct of the  $N_{\Delta_2}$ modulo $\Dscr$.
Let $A$ denote $P^*(N^*)$, $\abar^*$ denote the ultraproduct of the
$\abar^{\Delta_2}$, and $\langle \bbar_i,\bbar'_i\cbar_i:i\in I\rangle$ enumerate $Q^*(N^*)$.
Now we claim
\begin{claim}
\label{ultraproduct}
\begin{enumerate}
\item $\lg(\bbar_i) =\lg(\sigma(\bbar_i)) = \lg(\cbar)=n$ ; $\lg(\abar) = m$ and $\sigma^*(\bbar_i) =\bbar'_i$. 
\item $\langle \bbar_i\bbar'_i\cbar_i:i\in I\rangle$ is a sequence of $L(T)$-indiscernibles over
$\abar^*$.
\item For each finite $\Delta_2 \subseteq L(T)$ 
with $\Delta_2\subseteq F(\Delta_1,3n)$
and 
each finite subsequence  from  
 $\langle \bbar_i\bbar'_i\cbar_i:i\in I\rangle$ indexed by
$J$ 
of length
at least $k=f(\Delta_2,3n )$ the $\Delta_2$-type of
$\langle \bbar_i\bbar'_i\cbar_i:i\in J\rangle\abar$ is the $\Delta_2$-type of 
some  $(\sigma^{\Delta_2},\Delta_1, \Delta_2, \Delta_3,n)$-obstruction 
  $\langle\bbar^{\Delta_2},\sigma(\bbar^{\Delta_2}) \cbar^{\Delta_2}\rangle\abar^{\Delta_2}$
 in $M^{\Delta_2}$over the empty set.
\item $(\av_L(\langle \bbar'_i:i\in I\rangle/N) =\av_L(\langle\cbar_i:i\in I\rangle/N)$.
\end{enumerate}
\end{claim}
This claim follows directly from the properties of ultraproducts. 
 (For item 3, apply Fact~\ref{monot} and the definition of the ultrafilter
$D$.)

Let $\Gamma$ be the 
$L$-type in
 the variables $ \langle  \xbar_i \xbar'_i \ybar_i: i \in 
I\rangle\cup\{\zbar\}  $
over the empty set
of  $\langle \bbar_i,\bbar'_i\cbar_i:i\in I\rangle\abar^*$.
For any finite $\Delta \subset L(T)$, 
let $\chi_{\Delta,k}(\overline{\xbar},\overline{\xbar}'\overline{\ybar},
\zbar)$ be the $\Delta$-type  over  the empty set
 of a subsequence of
$k$ elements from $ \langle  \bbar_i \bbar'_i \cbar_i: i \in I\rangle$  
and $\abar^*$ from a realization of $\Gamma$ with $\zbar$ for $\abar^*$. 

Recall the definition of $\TP_{\Delta_1}$ from Notation~\ref{taudef}.
Let $r = f(\Delta_1,n)$ and let $\theta_{\Delta_1}( \xbar_0, \ldots 
\xbar_{r-1},
\xbar'_0, \ldots \xbar'_{r-1},
 \ybar_0, \ldots \ybar_{r-1},\zbar)$ be the 
 formula:\\
$(\exists \xbar,\xbar',\ybar) [
 \chi_{\Delta_1,r}
( \xbar_0, \ldots \xbar_{r-1},\xbar'_0, \ldots \xbar'_{r-1}, \ybar_0, 
\ldots \ybar_{r-1},\zbar)$\\
$
\wedge \TP_{\Delta_1}(\xbar,\xbar',\ybar, \xbar_0, \ldots \xbar_{r-1},\xbar'_0, \ldots \xbar'_{r-1}, \ybar_0, \ldots \ybar_{r-1},\zbar)
\wedge \sigma (\xbar'
) = \ybar]$.

Without loss of generality we assume $0, 1 \dots r-1$ index disjoint sequences.
Now we claim:

\begin{claim}
\label{ax}
If 
${\bf K_\sigma}$, the family of existentially closed models of $T_{\sigma}$,
 is axiomatized by $T^*_{\sigma}$,
$$T^*_\sigma \cup \Gamma 
\cup \{\sigma(\xbar_i) = \xbar'_i:i \in I\}
\vdash \theta_{\Delta_1}( \xbar_0, \ldots \xbar_{r-1},\xbar'_0, \ldots \xbar'_{r-1}, \ybar_0, \ldots \ybar_{r-1},\zbar).$$
\end{claim}
(Abusing notation we write this with the $\xbar_i,\xbar'_i,
\ybar_i$ free.)


For this, let 
$(M',\sigma') \models T^*_\sigma$
 such that 
 $\langle \bbar, \bbar'_i,\cbar_i:  i \in I\rangle\abar$ satisfy $\Gamma$ in $M'$
and  for each $i \in I$, $\sigma(\bbar_i) = \bbar'_i$.
Suppose $M' \restrict L \prec M''$ and $M''$ is an $|M'|^+$-saturated model of
$T$.  In $M''$ we can find $\bbar \bbar',\cbar$ realizing the average of 
  $\langle \bbar_i \bbar'_i \cbar_i:  i\in I\rangle$
over $M'$. 
Then $$\sigma'(\tp(\bbar/M')) =\sigma'( \av (\langle  \bbar_i: i \in I\rangle/M'))= 
\av (\langle  \bbar'_i: i \in I\rangle/M')$$
$$=
\av (\langle  \cbar_i:  i  \in I\rangle/M') = (\tp(\cbar/M')$$

(The first and last equalities are  by the choice of $\bbar,\bbar',\cbar$;
the second holds since for each $i$, $\sigma'(\bbar_i) = \bbar'_i$, the third
follows from clause 4 in the description of the ultraproduct.)
Now since $M''$ is $|M'|^+$-saturated there 
 is an automorphism $\sigma^{''}$ of $N$  extending $\sigma'$ and taking $\bbar$ to $\cbar$.

As $(M',\sigma') \models T^*_{\sigma}$, it 
is existentially closed.  So  we can pull $\bbar,\bbar'\cbar$
down to $M'$. Thus,
$(M',\sigma') \models \theta_{\Delta_1} ( \bbar_0, \ldots \bbar_{r-1},
\bbar'_0, \ldots \bbar'_{r-1},
 \cbar_0, \ldots \cbar_{r-1},\abar)$.  But $(M',\sigma')$ was an arbitrary model of  
$T^*_\sigma \cup \Gamma 
\cup \{\sigma(\xbar_i) = \xbar'_i: i \in I\}$; so 

$$T^*_\sigma \cup \Gamma\cup \{\sigma(\xbar_i) = \xbar'_i:i \in I\} \vdash $$
$$\theta_{\Delta_1}( \xbar_0, \ldots \xbar_{r-1},\xbar'_0, \ldots \xbar'_{r-1}, \ybar_0, \ldots \ybar_{r-1},\zbar).$$

 By compactness, some finite subset $\Gamma_0$ of $\Gamma$ and a finite number of the specifications of $\sigma$ suffice;
let  $\Delta^*$ be the formulas mentioned in  $\Gamma_0$ along with those in
$F(\Delta_1, 3n)$  and $k$ the number of
$x_i,y_i$ appearing in  $\Gamma_0$ 
and let $\Delta_2 = F(\Delta^*,n)$.
Without loss of generality, $k \geq f(\Delta_1,3n)$.
Then,
$T^*_\sigma  \vdash$ 
$$(\forall \xbar_0 \ldots \xbar_{k-1},\xbar'_0, \ldots \xbar'_{k-1}, \ybar_0, \ldots \ybar_{k-1})[(\chi_{\Delta_2,k}( \xbar_0, \ldots \xbar_{k-1},\xbar'_0, \ldots \xbar'_{k-1},
 \ybar_0, \ldots \ybar_{k-1},\zbar)$$
$$\wedge \bigwedge_{i<k} \sigma(\xbar_i) = \xbar'_i)\rightarrow\theta_{\Delta_1}( \xbar_0, \ldots \xbar_{r-1},\xbar'_0, \ldots \xbar'_{r-1}, \ybar_0, \ldots \ybar_{r-1},\zbar)].$$

By item 3) in Claim~\ref{ultraproduct}, 
fix a $\Delta_2$ and  $\Delta_3 = \Delta_3(\Delta_2,m)$ containing  $\Delta_2$ and 
$\langle\bbar^{\Delta_2}_{i}{\bbar'}^{\Delta_2}_{i}\cbar_{i}^{\Delta_2}:i< k\rangle \abar^{\Delta_2}$
which form a $(\sigma,\Delta_1,\Delta_2,\Delta_3,n)$-obstruction over $\abar^{\Delta_2}$ and so that:
 $$M^{\Delta_2} \models \chi_{\Delta_2,k}(\langle \bbar_i^{\Delta_2}{\bbar'}_i^{\Delta_2}\cbar_i^{\Delta_2}:i< k\rangle \abar^{\Delta_2} ).$$
By the definition of an obstruction, $\sigma(\bbar_i^{\Delta_2}) 
={\bbar'}_i^{\Delta_2}$.  So by the choice of $\Gamma_0$,

$$(M^{\Delta_2},\sigma^{\Delta_2})
 \models \theta_{\Delta_1}(\langle\bbar_i^{\Delta_2}{\bbar'}_i^{\Delta_2}\cbar_i^{\Delta_2}:i< r\rangle \abar^{\Delta_2}).$$

Now, let  $\bbar,\bbar'\cbar\in M^{\Delta_2}$ with 
$\sigma^{\Delta_2}(\bbar') = \cbar$ witness this sentence.
Then 
$$(M^{\Delta_2},\sigma^{\Delta_2}) 
\models\TP_{\Delta_1}(\langle\bbar_i^{\Delta_2}{\bbar'}_i^{\Delta_2}\cbar_i^{\Delta_2}:i< k\rangle \abar^{\Delta_2})$$
 By the definition of
obstruction, $$\sigma^{\Delta_2}(\tp_{\Delta_3}(\bbar'/A^{\Delta_2}\cup \bigcup_{i<k}\langle \bbar_i^{\Delta_2}\sigma(\bbar_i^{\Delta_2})\cbar_i\rangle \cup \abar^{\Delta_2}))
\neq \tp_{\Delta_3}(\cbar/A^{\Delta_2}\cup \bigcup_{i<k}\langle \bbar_i^{\Delta_2}\sigma(\bbar_i^{\Delta_2})\cbar_i^{\Delta_2}\rangle \cup \abar^{\Delta_2}).$$
This contradicts
that $\sigma^{\Delta_2}$ is an automorphism
 and we finish.

Finally we have the main result.

\begin{theorem}  
\label{final}
If $T$ is a stable theory, $T_\sigma$ has a model companion if
and only if 
$T_\sigma$ 
admits no $\sigma$-obstructions.
\end{theorem}

Proof.  
We showed in
Lemma~\ref{obsneeded} that if $T_{\sigma}$ has a model companion then there
is no obstruction. If there is no obstruction, 
Lemma~\ref{fcpobs} implies $T$ does not have the finite cover property.
 By Lemma~\ref{maindeal}
for every formula $\psi(\xbar,\ybar,\zbar)$ 
there is an $L$-formula $\theta_\psi(\zbar)$ (write out condition 3
 of Lemma~\ref{maindeal})
which for any $(M,\sigma) \models T_\sigma$ 
holds of any $\abar$ in $M$ if and only if  there exists $(N,\sigma)$,
   $(M,\sigma) \subseteq (N,\sigma)\models T_\sigma$
and $$(N,\sigma)\models (\exists \xbar \ybar)\psi(\xbar,\ybar,\abar) 
\wedge \sigma(\xbar) =\ybar.$$
Thus, the class of existentially closed models of $T_\sigma$ is 
axiomatized by the sentences:
$(\forall \zbar)\theta_\psi(\zbar) \rightarrow (\exists \xbar
 \ybar)\psi(\xbar,\ybar,\abar) \wedge \sigma(\xbar) =\ybar.$
(We can restrict to formulas of the form $\psi(\xbar,\sigma(\xbar),\abar)$ 
by the
standard trick (\cite{KikyoPillay, ChatHrushovski}).

Kikyo and Pillay \cite{KikyoPillay} note that if a strongly minimal
theory has the definable mulplicity property then $T_\sigma$ has a model
companion.  In view of Theorem~\ref{final}, this implies that if $T$ has
the definable multiplicity property, then $T_\sigma$ 
admits no $\sigma$-obstructions. 
Kikyo and Pillay conjecture that for a strongly minimal set, the converse
holds: if $T_\sigma$ has a model
companion then $T$ has
the definable multiplicity property.  They prove this result if $T$ is a 
finite cover of a theory with the finite multiplicity property.  It
would follow in general from a positive answer to the following
question.

\begin{question}  If the strongly minimal theory $T$ 
 with finite rank 
does not have the
definable multiplicity property, must it omit obstructions?
\end{question}

Pillay has given   a direct proof that if a strongly
minimal $T$ has
the definable multiplicity property, then $T_\sigma$ 
admits no $\sigma$-obstructions. 
Pillay has provided an insightful reworking  of the ideas here
in a note which is available on his website \cite{Pillayobs}.

Here is a final question:

\begin{question}  Can $T_\sigma$ for an $\aleph_0$-categorical $T$ admit
obstructions?
\end{question}


\end{document}